\documentclass[12pt, reqno]{amsart}
\usepackage{amsmath}
\usepackage{amsthm}
\usepackage{amssymb}
\usepackage{enumerate}
\usepackage{graphicx}
\usepackage{graphicx}
\usepackage{amssymb,amsfonts}
\usepackage{amsmath}
\usepackage[colorlinks,linkcolor=blue,anchorcolor=blue,citecolor=blue]{hyperref}
\usepackage[numbers,sort&compress]{natbib}
\usepackage{marvosym}
\usepackage{epstopdf}
\usepackage{caption}
\usepackage{multicol}
\usepackage{multirow}
\usepackage{makebox}
\usepackage{subfigure}
\numberwithin{equation}{section}
\newtheorem{theo}{Theorem} %[section]

\newtheorem{lem}{Lemma}
\newtheorem{cor}{Corollary}

\newtheorem{prop}{Proposition}

\begin{document}

\title{Some Sharp Chernoff type inequalities}
\author[Zhou and Zeng]{Yuqi Zhou$^1$\and Chunna Zeng$^{2,3}$$^*$}

\address{1.School of Mathematical Sciences,
 Chongqing Normal University,
Chongqing 401331, People's Republic of China}
\email{zhouyuqi202212@163.com}

\address{2.School of Mathematical Sciences, Chongqing Normal University, Chongqing 401331, People's Republic of China}
\email{zengchn@163.com}

%\address{3.Institut f\"{u}r Diskrete Mathematik und Geometrie, Technische Universit\"{a}t Wien, Wiedner Hauptstrasse 8-10/1046, 1040 Wien, Austria}
%\email{zengchn@163.com}

\thanks{This work was supported in part by the Major Special Project of the National Natural Science Foundation of China (Grant No. 12141101), the Young Top-Talent program of Chongqing (Grant No. CQYC2021059145), Natural Science Foundation Project of Chongqing (Grant No. cstc2020jcyj-msxmX0609), Technology Research Foundation of Chongqing Educational committee(Grant No. KJZD-K202200509).}
\thanks{{*}Corresponding author: Chunna Zeng}

%%%%%%%%%%%%%%%%%%%%%%%%%%%%%%%%%%%%%%%%%%%%%%%%%%%%%%%%%%%%%%%%%%%%%%%%%%%%%%%%%%%%%%%%%%%%%%%%%%%%%%%%%%%%
%%%%%%%%%%%%%%%%%%%%%%%%%%%%%%%%%%%%%%%%%%%%%%%%%%%%%%%%%%%%%%%%%%%%%%%%%%%%%%%%%%%%%%%%%%%%%%%%%%%%%%%%%%%
\begin{abstract}
 Two sharp Chernoff type inequalities are obtained for star body in $\mathbb{R}^2$, one of which is an extension of the dual Chernoff-Ou-Pan inequality, and the other is the reverse Chernoff type inequality. Furthermore, we establish a generalized dual symmetric mixed Chernoff inequality for two planar star bodies. As a direct consequence, a new proof of the the dual symmetric mixed isoperimetric inequality is presented.
\end{abstract}

\maketitle

%%%%%%%%%%%%%%%%%%%%%%%%%%%%%%%%%%%%%%%%%%%%%%%%%%%%%%%%%%%%%%%%%%%%%%%%%%%%%%%%%%%%%%%%%%%%%%%%%%%%%%%%%%%
%%%%%%%%%%%%%%%%%%%%%%%%%%%%%%%%%%%%%%%%%%%%%%%%%%%%%%%%%%%%%%%%%%%%%%%%%%%%%%%%%%%%%%%%%%%%%%%%%%%%%%%%%%%
\section{Introduction }
\par The isoperimetric inequality is one of the most classical geometric inequality. The classical isoperimetric inequality states that $L^2-4\pi A\geq 0$ for a domain with the boundary composing of the simple curve of area $A$ and length $L$, the equality holds if and only if this domain is a disc. Although this fact was already known to the ancient Greeks, its rigorous mathematical proof was not established until the variational method based on calculus appeared in the 19th century. Since then various proofs, generalizations and applications of the isoperimetric inequality were taken up (see \cite{Chernoff,Chavel,Dergiades,Osserman,Fusco,Figalli,Zhu-Li-Zhou,Zhu-Xu,Pan-Xu,Zeng-Ma-Zhou-Chen}).
The first geometric proof of the isoperimetric inequality was verified by Steiner via symmetrization, which can be founded in his work \cite{Steiner}.
In \cite{Fusco}, the quantitative sharp form of the isoperimetric inequality was proved by Fusco, Maggi and Pratelli. And the result deduced to a positive answer to a conjecture by Hall (\cite{Hall}).
In 1969, by employing the method due to Hurwitz, Chernoff \cite{Chernoff} demonstrated an isoperimetric type inequality involving the width function $w(P,\theta)$ of a convex body $P$, that is
\begin{equation}\label{1.1}
\frac{1}{2}\int_{0}^{\frac{\pi}{2}} w(P,\theta) w\left(P,\theta+\frac{\pi}{2}\right) d\theta \geq A(P),
\end{equation}
where $A(P)$ is the area of $P$, the equality holds if and only if $P$ is a disc. As the first milestone of the Chernoff type inequality, (\ref{1.1}) is a far-reaching generalization of the
isoperimetric inequality. The Chernoff type inequality (\ref{1.1}) exposes the crucial property of  the width function, and it reveals the well-known conclusion that a convex body with $A(P)>\frac{\pi}{4}$ must have $w(P,\theta)> 1$ in some direction.

\par In 2010, Ou and Pan \cite{Ou-Pan} introduced the $k$-order width function $w_k(P,\theta)$ of a convex body $P$
\begin{equation*}
w_{k}(P,\theta)=h(P,\theta)+h\left(P,\theta+\frac{2\pi}{k}\right)+\cdots+h\left(P,\theta+\frac{2(k-1)\pi}{k}\right),~~~k\geq2,
\end{equation*}
where $h(P,\theta)$ is the support function of $P$. They
gave some properties of convex bodies with $w_k(P,\theta)$ being constant and proposed an open problem motivated by the Blaschke-Lebesgue theorem. Furthermore, the authors obtained a generalization of (\ref{1.1}) involving the $k$-order width function $w_{k}(P,\theta)$, which is called the Chernoff-Ou-Pan inequality
\begin{equation}\label{1.2}
\frac{1}{k}\int_{0}^{\frac{\pi}{k}} w_k(P,\theta) w_k\left(P,\theta+\frac{\pi}{k}\right) d\theta \geq A(P),
\end{equation}
with equality if and only if $P$ is a disc. Note that for taking the limit as $k\rightarrow \infty$ in (\ref{1.2}), one can lead to a new proof of the classical isoperimetric inequality.
\par Let $P$, $Q$ be two planar convex bodies with area $A(P)$, $A(Q)$, respectively. One can find some elegant proofs that lead to the mixed isoperimetric type inequalities (see \cite{Xu-Zhou-Zhu,Zeng-Zhou-Yue}). In 1977, with the excellent work of Lutwak \cite{Lutwak1}, an inequality is established for the mixed width-integrals analogous to the Fenchel-Aleksandrov inequality for the mixed volumes, i.e,
\begin{equation}\label{1.3}
\frac{1}{8}\int_{0}^{2\pi} w(P,\theta) w(Q,\theta) d\theta \geq \sqrt{A(P)A(Q)},
\end{equation}
where $w(P,\theta)$, $w(Q,\theta)$ are the width functions of $P$, $Q$, respectively; with equality if and only if $P$ and $Q$ are discs. Inspired by the work of Ou and Pan \cite{Ou-Pan} and Lutwak \cite{Lutwak1}, the Chernoff-Ou-Pan inequality (\ref{1.2}) is generalized to the case of two convex bodies by Zhang \cite{Zhang1} as
\begin{equation}\label{1.4}
\frac{1}{2k}\int_{0}^{\frac{2\pi}{k}} w_k(P,\theta) w_k\left(Q,\theta+\frac{\pi}{k}\right) d\theta \geq \sqrt{A(P)A(Q)},
\end{equation}
with equality if and only if $P$ and $Q$ are discs. When $P=Q$, $k=2$, the author improved the Chernoff type inequality (\ref{1.1}) by combining the area of the Wigner caustic associated to a convex domain. Another byproduct of (\ref{1.4}) is the strengthened version of the classical isopermetric inequality. The interested readers can refer to \cite{Zhang1} for more details.

\par Actually, the support function of convex body in the Brunn-Minkowski theory corresponds to the radial function of star body in the dual theory. In this way, many geometric inequalities for star body have expanded rapidly (see \cite{Lutwak2,Lutwak4}).

\par The above cited work, especially in \cite{Ou-Pan}, make it
apparent that the time is ripe for the next step in the extension of radial function of star body. Zhang and Yang \cite{Zhang-Yang} defined the $k$-order radial function $\rho_{k}(S,\theta)$ of a planar star body $S$ by
\begin{equation*}
\rho_{k}(S,\theta)=\rho(S,\theta)+\rho\left(S,\theta+\frac{2\pi}{k}\right)+\cdots+\rho\left(S,\theta+\frac{2(k-1)\pi}{k}\right),~~~k\geq2,
\end{equation*}
where $\rho(S,\theta)$ is the radial function of $S$. Just as, $\rho_{2}(S,\theta)$ is the length of the chord passing through the origin in the direction $\vec{\mu}=(\cos\theta, \sin\theta)$, the function $\rho_{k}(S,\theta)$ is a natural generalization of the chord of $S$ passing through the origin. In particular, they presented that $\frac{C^2}{k^2}\leq A(S) <\frac{\pi C^2}{k}$ for star bodies with $\rho_{k}(S,\theta)=C$. Therefore, the Chernoff-Ou-Pan inequality (\ref{1.2}) was successfully extended from a convex body to a star body, then
\begin{equation}\label{1.5}
\frac{1}{k}\int_{0}^{\frac{\pi}{k}} \rho_k(S,\theta) \rho_k\left(S,\theta+\frac{\pi}{k}\right) d\theta \leq A(S),
\end{equation}
with equality if and only if the radial function of $S$ is the following form
\begin{align*}
\rho(\theta)=\frac{1}{2}a_0+\sum_{n=1}^{\infty}(a_{2nk}\cos 2nk\theta+b_{2nk}\sin 2nk\theta).
\end{align*}
(\ref{1.5}) is called the dual Chernoff-Ou-Pan inequality, and it
gave another proof of the dual isoperimetric inequality.
\par See also Mao and Yang \cite{Mao-Yang} for the expansion of (\ref{1.5}), the Chernoff type inequality for two planar star bodies $S$ and $T$ has reached to
\begin{equation}\label{1.6}
\frac{1}{2k^2}\int_{0}^{2\pi} \rho_k(S,\theta) \rho_k(T,\theta) d\theta \leq \sqrt{A(S)A(T)},
\end{equation}
with equality if and only if the radial functions of $S$ and $T$ have the same form
\begin{align*}
\rho(\theta)=\frac{1}{2}a_0+\sum_{n=1}^{\infty}(a_{nk}\cos nk\theta+b_{nk}\sin nk\theta).
\end{align*}
As a direct application of (\ref{1.6}), the dual symmetric mixed isoperimetric inequality $\tilde{A}(S,B)\tilde{A}(T,B)\leq\pi \sqrt{A(S)A(T)}$ for two planar star bodies can be obtained.

\par Recently, Fang ang Yang \cite{Fang-Yang1} proved the enhanced form of the Chernoff-Ou-Pan inequality (\ref{1.2}) by using Fourier series, see \cite{Fang-Yang1} for more details. A natural question is whether the dual Chernoff-Ou-Pan inequality (\ref{1.5}) can be strengthened. Using this insight, one of the aims of this article is to put forward the sharp version of the dual Chernoff-Ou-Pan inequality (\ref{1.5}) in a way modeled after \cite{Fang-Yang1}.

\begin{theo}\label{theo1}
Let $S$ be a planar star body with area $A(S)$. Suppose that $\frac{n}{k}$ is even, $a_n=b_n=0$.
For $k\geq2$, $k\in\mathbb{Z}^{+}$ and $0\leq\lambda\leq\frac{k}{\pi}$, the inequality
\begin{align}\label{1.7}
\int_{0}^{\frac{\pi}{k}} \rho_{k}(S,\theta) \rho_{k}\left(S,\theta+\frac{\pi}{k}\right)d\theta \leq kA(S)+\lambda\left(\tilde{A}(S,B)^2-\pi A(S) \right)
\end{align}
holds, where $\tilde{A}(S,B)$ is the dual mixed area of $S$ and an unit ball $B$ in $\mathbb{R}^2$; Moreover,
\par(i) if $0\leq\lambda<\frac{k}{\pi}$, then the equality in (\ref{1.7}) holds if and only if $S$ is a disc;
\par(ii) if $\lambda=\frac{k}{\pi}$, then the equality in (\ref{1.7}) holds if and only if the radial function of $S$ is of the form
\begin{align*}
\rho(S,\theta)=\frac{1}{2}a_0+a_1\cos \theta+b_1\sin\theta.
\end{align*}
\end{theo}
\par  For a planar star body, (\ref{1.7}) is stronger than the dual Chernoff-Ou-Pan inequality (\ref{1.5}) (see \cite{Zhang-Yang}), which can be seen from the classical dual isoperimetric inequality (\ref{2.7}). They further proposed the reverse Chernoff type inequality for a convex body  by using the reverse isoperimetric inequality (see \cite{Pan-Zhang}). As a consequence, it becomes natural to introduce the reverse isoperimetric inequality for a planar star body (see \cite{Fang}). We will prove the following theorem.
\begin{theo}\label{theo2}
Let $S$ be a planar star body with area $A(S)$. For $k\geq2$, $k\in\mathbb{Z}^{+}$ and $\mu\leq -k$, the inequality
\begin{align}\label{1.8}
\int_{0}^{\frac{\pi}{k}} \rho_{k}(S,\theta) \rho_{k}\left(S,\theta+\frac{\pi}{k}\right)d\theta \geq kA(S)+ \mu\tilde{A}(S)
\end{align}
holds, where $\tilde{A}(S)$ is the oriented area of $S$; Moreover,
\par (i) if $\mu< -k$, then the equality in (\ref{1.8}) holds if and only if $S$ is a disc;
\par (ii) if $\mu=-k$, then the equality in (\ref{1.8}) holds if and only if the radial function of $S$ is of the form
\begin{align*}
\rho(S,\theta)=\frac{1}{2}a_0+a_1\cos \theta+b_1\sin\theta.
\end{align*}
\end{theo}
\par Observe that (\ref{1.8}) is the reverse Chernoff type inequality involving the oriented area of a planar star body. In addition, we will establish the stabilities of (\ref{1.7}) and (\ref{1.8}) with the dual $L_2$ metrics in Section 3.

\par Another purpose of this paper is to generalize the dual symmetric mixed Chernoff inequality (\ref{1.6}) via the extension of independent variable.
\begin{theo}\label{theo3}
Let $S$, $T$ be two planar star bodies with area $A(S)$, $A(T)$, respectively. For $k\geq2$ and $k\in\mathbb{Z}^{+}$, then
\begin{align}\label{1.9}
\frac{1}{2k^2}\int_{0}^{2\pi} \rho_{k}(S,\theta) \rho_{k}(T,\theta+\alpha)d\theta \leq \sqrt{A(S)A(T)},~~~\alpha\in(0,2\pi).
\end{align}
The equality holds if and only if the radial functions of $S$ and $T$ have the same form
\begin{align*}
\rho(\theta)=\frac{1}{2}a_0+\sum_{n=1}^{\infty}(a_{nk}\cos nk\theta+b_{nk}\sin nk\theta).
\end{align*}
\end{theo}

\par One interesting feature of theorem \ref{theo3} is that it shows a novel proof of the dual symmetric mixed isoperimetric inequality $\tilde{A}(S,B)\tilde{A}(T,B)\leq\pi \sqrt{A(S)A(T)}$ (see \cite{Mao-Yang}).

\section{Preliminaries}

Denote by $\mathbb{R}^n$ the $n$-dimensional Euclidean space. Let $S$ be a compact star shaped set with respect to the origin, its radial function $\rho(S,\mu)$ is defined by
\begin{equation*}
\rho(S,\mu)=\max\left\{\lambda>0:\lambda\mu\in S\right\}
\end{equation*}
for $\mu\in\mathbb{S}^{n-1}$. $S$ is called a star body if $\rho(S,\mu)$ is continuous and positive.
\par Let $S$ be a planar star body. For $\theta\in[0,2\pi]$, define the area of $S$ and the oriented area of $S$ by
\begin{equation}\label{2.1}
A(S)=\frac{1}{2}\int_{0}^{2\pi}\rho^2(S,\theta)d\theta,
\end{equation}
\begin{equation}\label{2.2}
\tilde{A}(S)=\frac{1}{2}\int_{0}^{2\pi} \rho'^2(S,\theta) d\theta.
\end{equation}
It is clear that $\rho(S,\theta)$ is always continuous, bounded and $2\pi$-periodic. Thus, $\rho(S,\theta)$ can be expressed by Fourier series
\begin{equation}\label{2.3}
\rho(S,\theta)=\frac{1}{2}a_0+\sum_{n=1}^{\infty}(a_n\cos n\theta+b_n\sin n\theta),
\end{equation}
where
\begin{equation*}
a_0=\frac{1}{\pi}\int_{0}^{2\pi}\rho(S,\theta) d\theta
\end{equation*}
and
\begin{equation*}
a_n=\frac{1}{\pi}\int_{0}^{2\pi}\rho(S,\theta)\cos n\theta d\theta,~~~b_n=\frac{1}{\pi}\int_{0}^{2\pi}\rho(S,\theta)\sin n\theta d\theta,~~~n\in\mathbb{Z}^+.
\end{equation*}
By (\ref{2.1}), (\ref{2.2}), (\ref{2.3}) and the Parseval identity, we obtain
\begin{equation}\label{2.4}
A(S)=\frac{1}{4}\pi a_{0}^2+\frac{1}{2}\pi\sum_{n=1}^{\infty}(a_{n}^2+b_{n}^2),
\end{equation}
\begin{equation}\label{2.5}
\tilde{A}(S)=\frac{\pi}{2}\sum_{n=1}^{\infty} n^2(a_n ^2+b_n ^2).
\end{equation}

For a star body $S$ and an unit ball $B$ in $\mathbb{R}^2$, the dual mixed area $\tilde{A}(S,B)$ of $S$, $B$ is defined by
\begin{equation}\label{2.6}
\tilde{A}(S,B)=\frac{1}{2}\int_{0}^{2\pi}\rho(S,\theta)d\theta.
\end{equation}
Lutwak \cite{Lutwak4} proved the classical dual isoperimetric inequality which states that
\begin{equation}\label{2.7}
\tilde{A}(S,B)^2\leq \pi A(S),
\end{equation}
with equality if and only if $S$ is a disc. As a generalization of (\ref{2.7}), Mao and Yang \cite{Mao-Yang} obtained the dual symmetric mixed isoperimetric inequality
\begin{equation}\label{1}
\tilde{A}(S,B)\tilde{A}(T,B)\leq\pi \sqrt{A(S)A(T)},
\end{equation}
where $S$, $T$ are two planar star bodies.

\par In order to build stabilities of (\ref{1.7}) and (\ref{1.8}), the dual $L_2$ metric is needed
\begin{equation}\label{2.8}
\tilde{\delta_2}(S,T)=\left(\int_{0}^{2\pi}|\rho(S,\theta)-\rho(T,\theta)|^2d\theta\right)^{\frac{1}{2}}.
\end{equation}
\par In preparation for the proof of the main theorem, we shall use the following facts:
\begin{lem}\label{lem1}(\cite{Zhang-Yang})
Let $\rho_{k}(\theta)$ be the $k$-order radial function of a planar star body. For $k\geq2$, $k\in\mathbb{Z}^+$, then
\begin{equation}
\int_{0}^{\frac{\pi}{k}} \rho_{k}(\theta) \rho_{k}\left(\theta+\frac{\pi}{k}\right)d\theta =\frac{1}{2}\sum_{m=1}^{k}\int_{0}^{2\pi} \rho(\theta) \rho\left(\theta+\frac{2m-1}{k}\pi\right)d\theta.
\end{equation}
\end{lem}

\begin{lem}\label{lem2}
Let $f(\theta)$ and $g(\theta)$ be continuous, bounded and $2\pi$-periodic functions. For $k\geq2$, $k\in\mathbb{Z}^+$ and $\alpha\in(0,2\pi)$, then
\begin{align}
&\int_{0}^{2\pi} f_{k}(\theta) g_{k}(\theta+\alpha)d\theta \notag \\
=&\frac{k}{2}\sum_{m=1}^{k}\int_{0}^{2\pi} \left(f\left(\theta+\frac{2m\pi}{k}\right)g(\theta+\alpha)+f(\theta)g\left(\theta+\alpha+\frac{2m\pi}{k}\right)\right) d\theta,
\end{align}
where
\begin{equation*}
f_{k}(\theta)=f(\theta)+f\left(\theta+\frac{2\pi}{k}\right)+\cdots+f\left(\theta+\frac{2(k-1)\pi}{k}\right),
\end{equation*}
\begin{equation*}
g_{k}(\theta)=g(\theta)+g\left(\theta+\frac{2\pi}{k}\right)+\cdots+g\left(\theta+\frac{2(k-1)\pi}{k}\right).
\end{equation*}
\end{lem}
\emph{Proof.}
Since
\begin{align}\label{2}
&\int_{0}^{2\pi} f_{k}(\theta) g_{k}(\theta+\alpha)d\theta  \notag \\
=&\int_{0}^{2\pi} \sum_{i=1}^{k} f\left(\theta+\frac{2(i-1)\pi}{k}\right)\sum_{j=1}^{k} g\left(\theta+\alpha+\frac{2(j-1)\pi}{k}\right) d\theta \notag \\
=& \sum_{i, j=1}^{k} \int_{0}^{2\pi} f\left(\theta+\frac{2(i-1)\pi}{k}\right) g\left(\theta+\alpha+\frac{2(j-1)\pi}{k}\right) d\theta.
\end{align}
Let $\delta=\theta+\frac{2(i-1)\pi}{k}$, which together with the fact that $f(\theta)$ and $g(\theta)$ are $2\pi$-periodic functions yields
\begin{align*}
\int_{0}^{2\pi} f_{k}(\theta) g_{k}(\theta+\alpha)d\theta = \sum_{i, j=1}^{k} \int_{0}^{2\pi} f(\delta) g\left(\delta+\alpha+\frac{2(j-i)\pi}{k}\right) d\delta.
\end{align*}
Let $m=j-i$, we have
\begin{align}\label{3}
\int_{0}^{2\pi} f_{k}(\theta) g_{k}(\theta+\alpha)d\theta
=&k\sum_{m=0}^{k-1} \int_{0}^{2\pi} f(\delta) g\left(\delta+\alpha+\frac{2m\pi}{k}\right) d\delta \notag \\
=& k\sum_{m=1}^{k} \int_{0}^{2\pi} f(\theta) g\left(\theta+\alpha+\frac{2m\pi}{k}\right) d\theta.
\end{align}
Similarly, let $\delta=\theta+\frac{2(j-1)\pi}{k}$, it obtains
\begin{align}\label{4}
\int_{0}^{2\pi} f_{k}(\theta) g_{k}(\theta+\alpha)d\theta
=&\sum_{i, j=1}^{k} \int_{0}^{2\pi} f\left(\delta+\frac{2(i-j)\pi}{k}\right) g\left(\delta+\alpha\right) d\delta \notag \\
=& k \sum_{m=0}^{k-1} \int_{0}^{2\pi} f\left(\delta+\frac{2m\pi}{k}\right) g\left(\delta+\alpha\right) d\delta \notag \\
=& k \sum_{m=1}^{k} \int_{0}^{2\pi} f\left(\theta+\frac{2m\pi}{k}\right) g\left(\theta+\alpha\right) d\theta.
\end{align}
It follows from (\ref{2}), (\ref{3}) and (\ref{4}) that
\begin{align*}
&\int_{0}^{2\pi} f_{k}(\theta) g_{k}(\theta+\alpha)d\theta  \\
=&\frac{k}{2}\sum_{m=1}^{k}\int_{0}^{2\pi} \left(f\left(\theta+\frac{2m\pi}{k}\right)g(\theta+\alpha)+f(\theta)g\left(\theta+\alpha+\frac{2m\pi}{k}\right)\right) d\theta.
\end{align*}
\qed

\section{dual Chernoff-Ou-Pan inequalities and their stability properties}
\emph{Proof of Theorem 1.}\ \ \
By lemma \ref{lem1}, for $\rho_{k}(S,\theta)$ and $\rho_{k}\left(S,\theta+\frac{\pi}{k}\right)$, we obtain
\begin{align*}
\int_{0}^{\frac{\pi}{k}} \rho_{k}(S,\theta) \rho_{k}\left(S,\theta+\frac{\pi}{k}\right)d\theta =\frac{1}{2}\sum_{m=1}^{k}\int_{0}^{2\pi} \rho(S,\theta) \rho\left(S,\theta+\frac{2m-1}{k}\pi\right)d\theta,
\end{align*}
which together with (\ref{2.3}) yields
\begin{align*}
&\int_{0}^{2\pi} \rho(S,\theta) \rho\left(S,\theta+\frac{2m-1}{k}\pi\right)d\theta \\
=&\int_{0}^{2\pi} \left[ \frac{1}{2}a_0+\sum_{n=1}^{\infty}\left( a_n\cos n\theta+b_n\sin n\theta\right)\right]\\
&\left[ \frac{1}{2}a_0+\sum_{n=1}^{\infty}\left( a_n\cos n(\theta+\frac{2m-1}{k}\pi) +b_n\sin n(\theta+\frac{2m-1}{k}\pi)\right)\right]d\theta\\
=&\frac{\pi a_0^2}{2}+\pi \sum_{n=1}^{\infty}(a_n^2+b_n^2)\cos\frac{(2m-1)n\pi}{k}.
\end{align*}
Hence,
\begin{align*}
\int_{0}^{\frac{\pi}{k}} \rho_{k}(S,\theta) \rho_{k}\left(S,\theta+\frac{\pi}{k}\right)d\theta =\frac{k\pi a_0^2}{4}+\frac{\pi}{2} \sum_{n=1}^{\infty}(a_n^2+b_n^2)\sum_{m=1}^{k}\cos\frac{(2m-1)n\pi}{k}.
\end{align*}
When $n\neq kl,l\in\mathbb{Z}^{+}$, then
\begin{align*}
&\sum_{m=1}^{k}\cos\frac{(2m-1)n\pi}{k}\\
=&\frac{1}{\sin\frac{n\pi}{k}}\sum_{m=1}^{k}\left(\cos\frac{(2m-1)n\pi}{k}\sin\frac{n\pi}{k}\right)\\
=&\frac{1}{\sin\frac{n\pi}{k}}\sum_{m=1}^{k}\frac{1}{2}\left(\sin\frac{2mn\pi}{k}-\sin\frac{2(m-1)n\pi}{k}\right)\\
=&\frac{1}{\sin\frac{n\pi}{k}}\cdot \frac{1}{2}\left[(\sin\frac{2n\pi}{k}-0)+(\sin\frac{4n\pi}{k}-\sin\frac{2n\pi}{k})+\cdots+(\sin 2n\pi-\sin\frac{2(k-1)n\pi}{k})\right]\\
=&0,
\end{align*}
thus
\begin{align}\label{3.1}
\int_{0}^{\frac{\pi}{k}} \rho_{k}(S,\theta) \rho_{k}\left(S,\theta+\frac{\pi}{k}\right)d\theta =\frac{k\pi a_0^2}{4}+\frac{k\pi}{2} \sum_{l=1}^{\infty}(-1)^l(a_{kl}^2+b_{kl}^2).
\end{align}
From (\ref{2.4}) and (\ref{2.6}) it follows that
\begin{align}\label{3.2}
&\int_{0}^{\frac{\pi}{k}} \rho_{k}(S,\theta) \rho_{k}\left(S,\theta+\frac{\pi}{k}\right)d\theta - kA(S)-\lambda\left(\tilde{A}(S,B)^2-\pi A(S) \right)\notag \\
=&\frac{k\pi a_0^2}{4}+\frac{k\pi}{2} \sum_{l=1}^{\infty}(-1)^l(a_{kl}^2+b_{kl}^2)-\frac{\pi^2 a_0^2 \lambda}{4}+(\pi\lambda-k)\left(\frac{1}{4}\pi a_0^2+\frac{\pi}{2}\sum_{n=1}^{\infty}(a_n^2+b_n^2)\right)\notag \\
=&\frac{k\pi}{2} \sum_{l=1}^{\infty}(-1)^l(a_{kl}^2+b_{kl}^2)+\frac{k\pi}{2}(\frac{\pi\lambda}{k}-1)\sum_{n=1}^{\infty}(a_n^2+b_n^2)\notag \\
=&\frac{k\pi}{2}(\frac{\pi\lambda}{k}-1)(a_1^2+b_1^2)+\frac{k\pi}{2}\sum_{n=2}^{\infty}\left((\frac{\pi\lambda}{k}-1)+(-1)^{\frac{n}{k}}\right)(a_n^2+b_n^2).
\end{align}
Set $\varphi(\lambda)=\int_{0}^{\frac{\pi}{k}} \rho_{k}(S,\theta) \rho_{k}\left(S,\theta+\frac{\pi}{k}\right)d\theta - kA(S)-\lambda\left(\tilde{A}(S,B)^2-\pi A(S) \right)$. Observe that $\varphi(\lambda)$ is a monotone increasing function with respect to $\lambda$. When $\frac{n}{k}$ is even, $a_n=b_n=0$, then
\begin{align*}
\varphi\left(\frac{k}{\pi}\right)=\frac{k\pi}{2}\sum_{n=2}^{\infty}(-1)^{\frac{n}{k}}(a_n^2+b_n^2)\leq 0.
\end{align*}
Thus for any $0\leq\lambda\leq\frac{k}{\pi}$, we have $\varphi(\lambda)\leq \varphi\left(\frac{k}{\pi}\right)\leq 0$.\\
(i) If $0\leq\lambda<\frac{k}{\pi}$, then $\frac{\pi\lambda}{k}-1<0$ and $(\frac{\pi\lambda}{k}-1)+(-1)^{\frac{n}{k}}<0$. It is obviously that the equality in (\ref{1.7}) holds if and only if $a_n=b_n=0$ for all $n\geq1$, that is, $S$ is a disc.\\
(ii) If $\lambda=\frac{k}{\pi}$, then the equality in (\ref{1.7}) holds if and only if $a_n=b_n=0$ for all $n\geq2$. In other words, the radial function of $S$ is of the form $\rho(S,\theta)=\frac{1}{2}a_0+a_1\cos \theta+b_1\sin\theta$.
\qed

\emph{Remark 1.} It follows from (\ref{2.7}) that for $0\leq\lambda\leq\frac{k}{\pi}$, thus
\begin{align*}
 kA(S)+\lambda\left(\tilde{A}(S,B)^2-\pi A(S) \right)\leq kA(S).
\end{align*}
As a direct consequence of Theorem \ref{theo1}, we obtain

\begin{cor}\cite{Zhang-Yang}
Let $S$ be a planar star body with area $A(S)$, then
\begin{align}\label{3.3}
\frac{1}{k}\int_{0}^{\frac{\pi}{k}} \rho_{k}(S,\theta) \rho_{k}\left(S,\theta+\frac{\pi}{k}\right)d\theta \leq A(S),
\end{align}
with the equality in (\ref{3.3}) holds if and only if the radial function of $S$ is of the form
\begin{align*}
\rho(S,\theta)=\frac{1}{2}a_0+\sum_{n=1}^{\infty}(a_{2nk}\cos 2nk\theta+b_{2nk}\sin 2nk\theta).
\end{align*}
\end{cor}

\emph{Proof of Theorem 2.}\ \ \
By (\ref{2.4}),(\ref{2.5}) and (\ref{3.1}), we have
\begin{align}\label{3.4}
&\int_{0}^{\frac{\pi}{k}} \rho_{k}(S,\theta) \rho_{k}\left(S,\theta+\frac{\pi}{k}\right)d\theta -kA(S)- \mu\tilde{A}(S)\notag\\
=&\frac{k\pi a_0^2}{4}+\frac{k\pi}{2} \sum_{l=1}^{\infty}(-1)^l(a_{kl}^2+b_{kl}^2)-k\left(\frac{1}{4}\pi a_{0}^2+\frac{1}{2}\pi\sum_{n=1}^{\infty}(a_{n}^2+b_{n}^2)\right)\notag\\
&-\mu\left(\frac{\pi}{2}\sum_{n=1}^{\infty} n^2(a_n ^2+b_n ^2)\right)\notag\\
=&\frac{k\pi}{2} \sum_{n=2}^{\infty}(-1)^{\frac{n}{k}}(a_n^2+b_n^2)+\frac{k\pi}{2} \sum_{n=1}^{\infty}\left(-1-\frac{\mu n^2}{k}\right)(a_n^2+b_n^2)\notag\\
=&\frac{k\pi}{2} \sum_{n=2}^{\infty}\left[(-1)^{\frac{n}{k}}+\left(-1-\frac{\mu n^2}{k}\right)\right](a_n^2+b_n^2)+\frac{k\pi}{2}\left(-1-\frac{\mu}{k}\right)(a_1^2+b_1^2).
\end{align}
Set $\psi(\mu)=\int_{0}^{\frac{\pi}{k}} \rho_{k}(S,\theta) \rho_{k}\left(S,\theta+\frac{\pi}{k}\right)d\theta -kA(S)- \mu\tilde{A}(S)$. Observe that $\psi(\mu)$ is a monotone decreasing function with respect to $\mu$ such that
\begin{align*}
\psi(-k)=\frac{k\pi}{2}\sum_{n=2}^{\infty}\left[(-1)^{\frac{n}{k}}-1+n^2\right](a_n^2+b_n^2) \geq 0.
\end{align*}
Thus for any $\mu\leq-k$, one obtain $\psi(\mu)\geq \psi(-k)\geq0$.\\
(i) If $\mu<-k$, then $-1-\frac{\mu}{k}>0$ and $(-1)^{\frac{n}{k}}+\left(-1-\frac{\mu n^2}{k}\right)>0$. Naturally, the equality in (\ref{1.8}) holds if and only if $a_n=b_n=0$ for all $n\geq1$, that is, $S$ is a disc.\\
(ii) If $\mu=-k$, then the equality in (\ref{1.8}) holds if and only if $a_n=b_n=0$ for all $n\geq2$. In other words, the radial function of $S$ is of the form $\rho(S,\theta)=\frac{1}{2}a_0+a_1\cos \theta+b_1\sin\theta$.
\qed

\par Next, we use the dual $L_2$ metrics to build the stabilities of (\ref{1.7}) and (\ref{1.8}).

\begin{theo}
Let $S$ be a planar star body with area $A(S)$.
Suppose that $\frac{n}{k}$ is even, $a_n=b_n=0$.
For $k\geq2$, $k\in\mathbb{Z}^{+}$ and $0\leq\lambda\leq\frac{k}{\pi}$, then
\begin{align}\label{3.5}
&\int_{0}^{\frac{\pi}{k}} \rho_{k}(S,\theta) \rho_{k}\left(S,\theta+\frac{\pi}{k}\right)d\theta- kA(S)-\lambda\left(\tilde{A}(S,B)^2-\pi A(S) \right)\notag \\
&\leq \frac{1}{2}(\pi \lambda-k)\tilde{\delta_2}\left(S,\frac{a_0}{2} B\right)^2,
\end{align}
where $\tilde{A}(S,B)$ is the dual mixed area of $S$ and $B$; Moreover, the equality in (\ref{3.5}) holds if and only if the radial function of $S$ is of the form
\begin{align*}
\rho(S,\theta)=\frac{1}{2}a_0+a_1\cos \theta+b_1\sin\theta.
\end{align*}
\end{theo}

\emph{Proof.}
By (\ref{2.8}) and the Parserval equality, we have
\begin{align}\label{3.6}
\tilde{\delta_2}\left(S,\frac{a_0}{2} B\right)^2=&\int_{0}^{2\pi}|\rho(S,\theta)-\rho(\frac{a_0}{2} B,\theta)|^2d\theta\notag\\
=&\int_{0}^{2\pi}\left|\frac{1}{2}a_0+\sum_{n=1}^{\infty}(a_n\cos n\theta+b_n\sin n\theta)-\frac{1}{2}a_0\right|^2d\theta\notag\\
=&\int_{0}^{2\pi}\left|\sum_{n=1}^{\infty}(a_n\cos n\theta+b_n\sin n\theta)\right|^2d\theta\notag\\
=&\pi\sum_{n=1}^{\infty}(a_n^2+b_n^2).
\end{align}
Combing (\ref{3.2}) and (\ref{3.6}) yields that
\begin{align*}
&\int_{0}^{\frac{\pi}{k}} \rho_{k}(S,\theta) \rho_{k}\left(S,\theta+\frac{\pi}{k}\right)d\theta- kA(S)-\lambda\left(\tilde{A}(S,B)^2-\pi A(S) \right)\\
=&\frac{k\pi}{2}(\frac{\pi\lambda}{k}-1)(a_1^2+b_1^2)+\frac{k\pi}{2}\sum_{n=2}^{\infty}\left((\frac{\pi\lambda}{k}-1)+(-1)^{\frac{n}{k}}\right)(a_n^2+b_n^2)\\
=&\frac{k\pi}{2}(\frac{\pi\lambda}{k}-1)(a_1^2+b_1^2)+\frac{k\pi}{2}\sum_{n=2}^{\infty}\left((\frac{\pi\lambda}{k}-1)-1\right)(a_n^2+b_n^2)\\
\leq&\frac{k\pi}{2}(\frac{\pi\lambda}{k}-1)(a_1^2+b_1^2)+\frac{k\pi}{2}\sum_{n=2}^{\infty}\left(\frac{\pi\lambda}{k}-1\right)(a_n^2+b_n^2)\\
=&\frac{k\pi}{2}(\frac{\pi\lambda}{k}-1)\sum_{n=1}^{\infty}(a_n^2+b_n^2)\\
=&\frac{1}{2}(\pi \lambda-k)\tilde{\delta_2}\left(S,\frac{a_0}{2} B\right)^2,
\end{align*}
with equality if and only if $a_n=b_n=0$ for $n\geq 2$, that is, the radial function of $S$ is of the form $\rho(S,\theta)=\frac{1}{2}a_0+a_1\cos \theta+b_1\sin\theta$.
\qed

\begin{theo}
Let $S$ be a planar star body with area $A(S)$. For $k\geq2$, $k\in\mathbb{Z}^{+}$ and $\mu\leq -k$, then
\begin{align}\label{3.7}
\int_{0}^{\frac{\pi}{k}} \rho_{k}(S,\theta) \rho_{k}\left(S,\theta+\frac{\pi}{k}\right)d\theta - kA(S)- \mu\tilde{A}(S)\geq\left(-k-\frac{\mu}{2}\right)\tilde{\delta_2}\left(S,\frac{a_0}{2} B\right)^2,
\end{align}
where $\tilde{A}(S)$ is the oriented area of $S$; Moreover, the equality in (\ref{3.7}) holds if and only if $S$ is a disc.
\end{theo}

\emph{Proof.}
From (\ref{3.4}) and (\ref{3.6}), it follows that
\begin{align*}
&\int_{0}^{\frac{\pi}{k}} \rho_{k}(S,\theta) \rho_{k}\left(S,\theta+\frac{\pi}{k}\right)d\theta - kA(S)- \mu\tilde{A}(S)\\
=&\frac{k\pi}{2} \sum_{n=2}^{\infty}\left[(-1)^{\frac{n}{k}}+\left(-1-\frac{\mu n^2}{k}\right)\right](a_n^2+b_n^2)+\frac{k\pi}{2}\left(-1-\frac{\mu}{k}\right)(a_1^2+b_1^2)\\
\geq&\frac{k\pi}{2} \sum_{n=2}^{\infty}\left(-2-\frac{\mu n^2}{k}\right)(a_n^2+b_n^2)+\frac{k\pi}{2}\left(-1-\frac{\mu}{k}\right)(a_1^2+b_1^2)\\
\geq&\frac{k\pi}{2} \sum_{n=2}^{\infty}\left(-2-\frac{\mu}{k}\right)(a_n^2+b_n^2)+\frac{k\pi}{2}\left(-2-\frac{\mu}{k}\right)(a_1^2+b_1^2)\\
=&\frac{k\pi}{2} \left(-2-\frac{\mu}{k}\right)\sum_{n=1}^{\infty}(a_n^2+b_n^2)\\
=&\left(-k-\frac{\mu}{2}\right)\tilde{\delta_2}\left(S,\frac{a_0}{2} B\right)^2,
\end{align*}
with the equality in (\ref{3.7}) holds if and only if $a_n=b_n=0$ for $n\geq 1$, that is, $S$ is a disc.
\qed

\section{A generalized dual symmetric mixed Chernoff inequality}
In this section, we will deal with a generalized dual symmetric mixed Chernoff inequality of two star bodies.

\emph{Proof of Theorem 3.}\ \ \
From Lemma \ref{lem2} and the Schwartz inequality, one has
\begin{align*}
&\frac{1}{2k^2}\int_{0}^{2\pi} \rho_{k}(S,\theta) \rho_{k}(T,\theta+\alpha)d\theta\\
=&\frac{1}{4k}\sum_{m=1}^{k}\int_{0}^{2\pi} \left(\rho(S,\theta+\frac{2m\pi}{k})\rho(T,\theta+\alpha)+\rho(S,\theta)\rho(T,\theta+\alpha+\frac{2m\pi}{k})\right)d\theta\\
\leq&\frac{1}{2k}\sum_{m=1}^{k}\left(\frac{1}{2}\int_{0}^{2\pi}\rho^{2}(S,\theta+\frac{2m\pi}{k})d\theta\right)^{\frac{1}{2}}\cdot\left(\frac{1}{2}\int_{0}^{2\pi}\rho^{2}(T,\theta+\alpha)d\theta\right)^{\frac{1}{2}}\\
&+\frac{1}{2k}\sum_{m=1}^{k}\left(\frac{1}{2}\int_{0}^{2\pi}\rho^{2}(S,\theta)d\theta\right)^{\frac{1}{2}}\cdot\left(\frac{1}{2}\int_{0}^{2\pi}\rho^{2}(T,\theta+\alpha+\frac{2m\pi}{k})d\theta\right)^{\frac{1}{2}}\\
=&\sqrt{A(S)A(T)},
\end{align*}
where the equality holds if and only if for $m=1,2,\cdots,k$,
\begin{align*}
\rho\left(S,\theta+\frac{2m\pi}{k}\right)=c_1\rho(T,\theta+\alpha),
\end{align*}
\begin{align*}
\rho(S,\theta)=c_2\rho\left(T,\theta+\alpha+\frac{2m\pi}{k}\right),
\end{align*}
where $c_1$ and $c_2$ are two constants. It is easy to check that
\begin{align*}
\rho(S,\theta)&=c_2\rho\left(T,\theta+\alpha+\frac{2m\pi}{k}\right)=\frac{c_2}{c_1}\rho(S,\theta+\frac{2 \cdot 2m\pi}{k})=\frac{c_2^2}{c_1}\rho\left(T,\theta+\alpha+\frac{3\cdot2m\pi}{k}\right)\\
&=\cdots=\left(\frac{c_2}{c_1}\right)^{\frac{k}{2}}\rho(S,\theta+2m\pi)=\left(\frac{c_2}{c_1}\right)^{\frac{k}{2}}\rho(S,\theta),
\end{align*}
which together with the fact $\rho(S,\theta)>0$ yields $\frac{c_2}{c_1}=1$. Hence, the radial functions of $S$ and $T$ satisfy
\begin{align*}
\rho(\theta)=\rho\left(\theta+\frac{4m\pi}{k}\right).
\end{align*}
Moreover, (\ref{2.3}) leads that
\begin{align*}
\sum_{n=1}^{\infty}a_n\left(\cos n\theta-\cos n(\theta+\frac{4m\pi}{k})\right)+\sum_{n=1}^{\infty}b_n\left(\sin n\theta-\sin n(\theta+\frac{4m\pi}{k})\right)=0,
\end{align*}
by a direct computation, one obtain
\begin{align*}
\sum_{n=1}^{\infty} a_n\sin \frac{2mn\pi}{k} \sin \left(n\theta+\frac{2mn\pi}{k}\right)-\sum_{n=1}^{\infty} b_n\sin \frac{2mn\pi}{k} \cos \left(n\theta+\frac{2mn\pi}{k}\right)=0.
\end{align*}
Thus, $a_n\sin \frac{2mn\pi}{k}=b_n\sin \frac{2mn\pi}{k}=0$ for $n\in \mathbb{Z}^+$, which implies $a_n=b_n=0$ when $n\neq kl$, $l\in \mathbb{Z}^+$. That is, the radial functions of $S$ and $T$ have the same form
\begin{align*}
\rho(\theta)=\frac{1}{2}a_0+\sum_{n=1}^{\infty}(a_{nk}\cos nk\theta+b_{nk}\sin nk\theta).
\end{align*}
\qed

\begin{prop}\label{prop1}
Let $S$, $T$ be two planar star bodies. For $k\geq2$ and $k\in\mathbb{Z}^{+}$, then,
\begin{align}\label{4.3}
\lim_{k\rightarrow\infty}\frac{1}{2k^2}\int_{0}^{2\pi} \rho_{k}(S,\theta) \rho_{k}(T,\theta+\alpha)d\theta=\frac{1}{\pi} \tilde{A}(S,B)\tilde{A}(T,B),~~~\alpha\in(0,2\pi),
\end{align}
where $\tilde{A}(S,B)$ is the dual mixed area of $S$ and $T$ and $\tilde{A}(T,B)$ is that of $S$ and $T$.
\end{prop}

\emph{Proof.}\ \ \
For all $\xi\in(0,2\pi)$, there exists an $i\in[0,k-1]$ such that $\xi\in\left(\frac{2i\pi}{k},\frac{2(i+1)\pi}{k}\right)$. From the definition of definite integral and periodicity of $\rho(\theta)$, we have
\begin{align}\label{4.4}
\lim_{k\rightarrow\infty}\frac{2\pi}{k}\rho_{k}(S,\xi)&=\lim_{k\rightarrow\infty}\frac{2\pi}{k}\sum_{m=0}^{k-1}\rho\left(S,\xi+\frac{2m\pi}{k}\right)\notag\\
&=\lim_{k\rightarrow\infty}\frac{2\pi}{k}\left(\sum_{m=1}^{k}\rho\left(S,\xi+\frac{2m\pi}{k}\right)+\rho(S,\xi)-\rho(S,\xi+2\pi)\right)\notag\\
&=\lim_{k\rightarrow\infty}\frac{2\pi}{k}\sum_{m=1}^{k}\rho\left(S,\xi+\frac{2m\pi}{k}\right)\notag\\
&=\int_{0}^{2\pi}\rho(S,\theta)d\theta.
\end{align}
Similarly,
\begin{align}\label{4.5}
\lim_{k\rightarrow\infty}\frac{2\pi}{k}\rho_{k}(T,\xi+\alpha)=\int_{0}^{2\pi}\rho(T,\theta)d\theta.
\end{align}
Hence, using (\ref{4.4}), (\ref{4.5}) and the Cauchy mean-value theorem, there is a $\zeta\in(0,2\pi)$ such that
\begin{align*}
&\lim_{k\rightarrow\infty}\frac{1}{2k^2}\int_{0}^{2\pi} \rho_{k}(S,\theta) \rho_{k}(T,\theta+\alpha)d\theta\\
=&\lim_{k\rightarrow\infty}\frac{1}{2k^2}\cdot 2\pi  \rho_{k}(S,\zeta) \rho_{k}(T,\zeta+\alpha)\\
=&\frac{1}{4\pi}\left(\lim_{k\rightarrow\infty}\frac{2\pi}{k}\rho_{k}(S,\zeta)\right)\cdot\left(\lim_{k\rightarrow\infty}\frac{2\pi}{k}\rho_{k}(T,\zeta+\alpha)\right)\\
=&\frac{1}{\pi}\tilde{A}(S,B)\tilde{A}(T,B).
\end{align*}
\qed

\emph{Remark 2.} Combining Theorem \ref{theo3} with Proposition \ref{prop1} shows another proof of the the dual symmetric mixed isoperimetric inequality $\tilde{A}(S,B)\tilde{A}(T,B)\leq\pi \sqrt{A(S)A(T)}$.

\begin{prop}
Let $S$, $T$ be two planar star bodies with area $A(S)$, $A(T)$, respectively, then
\begin{align}\label{4.6}
\tilde{A}(S,B)\tilde{A}(T,B)\leq\pi \sqrt{A(S)A(T)},
\end{align}
where $\tilde{A}(S,B)$ is the dual mixed area of $S$ and $B$ and $\tilde{A}(T,B)$ is that of $T$ and $B$; Moreover, the equality in (\ref{4.6}) holds if and only if the radial functions of $S$ and $T$ have the same form
\begin{align*}
\rho(\theta)=\frac{1}{2}a_0+\sum_{n=1}^{\infty}(a_{nk}\cos nk\theta+b_{nk}\sin nk\theta).
\end{align*}
\end{prop}

\emph{Proof.}\ \ \
It follows from (\ref{1.9}) that for $k\geq2$ and $k\in\mathbb{Z}^{+}$, we have
\begin{align*}
\frac{1}{2k^2}\int_{0}^{2\pi} \rho_{k}(S,\theta) \rho_{k}(T,\theta+\alpha)d\theta \leq \sqrt{A(S)A(T)},
\end{align*}
which together with (\ref{4.3}) yields that
\begin{align*}
\lim_{k\rightarrow\infty}\frac{1}{2k^2}\int_{0}^{2\pi} \rho_{k}(S,\theta) \rho_{k}(T,\theta+\alpha)d\theta=\frac{1}{\pi} \tilde{A}(S,B)\tilde{A}(T,B).
\end{align*}
Thus,
\begin{align*}
\tilde{A}(S,B)\tilde{A}(T,B)\leq\pi \sqrt{A(S)A(T)},
\end{align*}
with the equality in (\ref{4.6}) holds if and only if the radial functions of $S$ and $T$ have the same form $\rho(\theta)=\frac{1}{2}a_0+\sum_{n=1}^{\infty}(a_{nk}\cos nk\theta+b_{nk}\sin nk\theta)$.
\qed

%%%%%%%%%%%%%%%%%%%%%%%%%%%%%%%%%%%%%%%%%%%%%%%%%%%%%%%%%%%%%%%%%%%%%%%%%%%%%%%%%%%%%%%%%%%%%%%%%%%%%%%%%%%

%%%%%%%%%%%%%%%%%%%%%%%%%%%%%%%%%%%%%%%%%%%%%%%%%%%%%%%%%%%%%%%%%%%%%%%%%%%%%%%%%%%%%%%%%%%%%%%%%%%%%%%%%%%

%%%%%%%%%%%%%%%%%%%%%%%%%%%%%%%%%%%%%%%%%%%%%%%%%%%%%%%%%%%%%%%%%%%%%%%%%%%%%%%%%%%%%%%%%%%%%%%%%%%%%%%%%%%
%%%%%%%%%%%%%%%%%%%%%%%%%%%%%%%%%%%%%%%%%%%%%%%%%%%%%%%%%%%%%%%%%%%%%%%%%%%%%%%%%%%%%%%%%%%%%%%%%%%%%%%%%%%

\end{document}